\documentclass[12pt]{amsart} 

\newcommand{\eeproof}{\hspace{10pt} \ \hfill $\Box$
\medskip
}
\usepackage{amsfonts}
\usepackage{latexsym}
\usepackage{pstricks}

\catcode`\@=11
\long\def\@makecaption#1#2{%
    \vskip 10pt

\setbox\@tempboxa\hbox{
      \small\sf{\bfcaptionfont #1. }\ignorespaces #2}%
    \ifdim \wd\@tempboxa >\captionwidth {%
        \rightskip=\@captionmargin\leftskip=\@captionmargin
        \unhbox\@tempboxa\par}%
      \else
        \hbox to\hsize{\hfil\box\@tempboxa\hfil}%
    \fi
}
\font\bfcaptionfont=cmssbx10 scaled \magstephalf
\newdimen\@captionmargin\@captionmargin=2\parindent
\newdimen\captionwidth\captionwidth=\hsize
\catcode`\@=12

\newcommand{\MILN}{
\psset{unit=.5cm}
\begin{pspicture}[.4](-.3,-.3)(1.5,1.3)
\psline[linewidth=.5pt,linestyle=dotted,dotsep=1pt](.1,.1)(.6,.5)
\psline[linewidth=.5pt,linestyle=dotted,dotsep=1pt](.1,.9)(.6,.5)
\psline[linewidth=.5pt,linestyle=dotted,dotsep=1pt](.6,.5)(1.1,.5)
\end{pspicture}
}

\newcommand{\MUa}{
\begin{pspicture}[.4](-.5,-.5)(1.5,1.5)
\psline[linewidth=.5pt,linestyle=dotted,dotsep=1pt](0,0)(0,1)
\psecurve[linewidth=.5pt,linestyle=dotted,dotsep=1pt](-.5,.5) (0,.5)
(.3,.5)(1,0) (1.1,-.5) 

\put(-.4,.9){$\scriptstyle{i}$}
\put(-.4,0){$\scriptstyle{j}$}
\put(1.2,.9){$\scriptstyle{l}$}
\put(1.2,0){$\scriptstyle{k}$}
\psdots[dotscale=.6](0,0)(0,1)(1,0)(1,1)
\end{pspicture}
}

\newcommand{\MUb}{
\begin{pspicture}[.4](-.5,-.5)(1.5,1.5)
\psline[linewidth=.5pt,linestyle=dotted,dotsep=1pt](0,0)(0,1)
\psecurve[linewidth=.5pt,linestyle=dotted,dotsep=1pt](-.5,.5) (0,.5)
(.3,.5) (1,1) (1.1,1.5) 

\put(-.4,.9){$\scriptstyle{i}$}
\put(-.4,0){$\scriptstyle{j}$}
\put(1.2,.9){$\scriptstyle{l}$}
\put(1.2,0){$\scriptstyle{k}$}
\psdots[dotscale=.6](0,0)(0,1)(1,0)(1,1)
\end{pspicture}
}

\newcommand{\MUc}{
\begin{pspicture}[.4](-.5,-.5)(1.5,1.5)
\psline[linewidth=.5pt,linestyle=dotted,dotsep=1pt](1,0)(1,1)
\psecurve[linewidth=.5pt,linestyle=dotted,dotsep=1pt](-.1,-.5)(0,0)
(.7,.5)(1,.5)(1.5,.5) 

\put(-.4,.9){$\scriptstyle{i}$}
\put(-.4,0){$\scriptstyle{j}$}
\put(1.2,.9){$\scriptstyle{l}$}
\put(1.2,0){$\scriptstyle{k}$}
\psdots[dotscale=.6](0,0)(0,1)(1,0)(1,1)
\end{pspicture}
}

\newcommand{\MUd}{
\begin{pspicture}[.4](-.5,-.5)(1.5,1.5)
\psline[linewidth=.5pt,linestyle=dotted,dotsep=1pt](1,0)(1,1)
\psecurve[linewidth=.5pt,linestyle=dotted,dotsep=1pt](-.1,1.5)(0,1)
(.7,.5)(1,.5)(1.5,.5) 
\put(-.4,.9){$\scriptstyle{i}$}
\put(-.4,0){$\scriptstyle{j}$}
\put(1.2,.9){$\scriptstyle{l}$}
\put(1.2,0){$\scriptstyle{k}$}
\psdots[dotscale=.6](0,0)(0,1)(1,0)(1,1)
\end{pspicture}
}

\newcommand{\FIGEX}{
\begin{pspicture}[.4](-.5,-.5)(4.5,3.5)
\pscurve[linewidth=.5pt,arrows=->,arrowscale=2](.4,1.45)(1,1.5)
        (1.1,1)(.7,0)(1,-1)(2,-.5)(3,-1)(3.5,.5)(2.5,1)(3.5,1)(4.5,1)
        (4,3.5)(2.5,3)(0,4)(-.5,2)  
\psline[linewidth=.5pt,linestyle=dotted,dotsep=1pt](1,0)(2,.5)
\psline[linewidth=.5pt,linestyle=dotted,dotsep=1pt](2,.5)(3,0)
\psline[linewidth=.5pt,linestyle=dotted,dotsep=1pt](2,.5)(2,1.5)
\psline[linewidth=.5pt,linestyle=dotted,dotsep=1pt](2,1.5)(3,2)
\psline[linewidth=.5pt,linestyle=dotted,dotsep=1pt](2,1.5)(1,2)
\psline[linewidth=.5pt,linestyle=dotted,dotsep=1pt](3,2)(3.5,3)
\psline[linewidth=.5pt,linestyle=dotted,dotsep=1pt](3,2)(4,1.5)
\psline[linewidth=.5pt,linestyle=dotted,dotsep=1pt](1,2)(.5,3)
\psline[linewidth=.5pt,linestyle=dotted,dotsep=1pt](1,2)(0,1.5)
\put(-.4,1.3){$\scriptstyle{1}$}
\put(.1,2.9){$\scriptstyle{5}$}
\put(1.5,1.1){$\scriptstyle{4}$}
\put(.9,-.6){$\scriptstyle{3}$}
\put(2.9,-.6){$\scriptstyle{2}$}
\put(4.1,1.3){$\scriptstyle{7}$}
\put(3.6,2.9){$\scriptstyle{6}$}
\psdots[dotscale=.6](0,1.5)(.5,3)(2,1.5)(3.5,3)(4,1.5)(1,0)(3,0)
\end{pspicture}
}

\newcommand{\FIGEXi}{
\begin{pspicture}[.4](-.5,1)(4.5,3.5)
\psline[linewidth=.5pt,linestyle=dotted,dotsep=1pt](2,1.5)(3,2)
\psline[linewidth=.5pt,linestyle=dotted,dotsep=1pt](2,1.5)(1,2)
\psline[linewidth=.5pt,linestyle=dotted,dotsep=1pt](3,2)(3.5,3)
\psline[linewidth=.5pt,linestyle=dotted,dotsep=1pt](3,2)(4,1.5)
\psline[linewidth=.5pt,linestyle=dotted,dotsep=1pt](1,2)(.5,3)
\psline[linewidth=.5pt,linestyle=dotted,dotsep=1pt](1,2)(0,1.5)
\put(-.4,1.3){$\scriptstyle{3}$}
\put(.1,2.9){$\scriptstyle{2}$}
\put(1.8,.9){$\scriptstyle{1}$}
\put(4.1,1.3){$\scriptstyle{4}$}
\put(3.6,2.9){$\scriptstyle{5}$}
\psdots[dotscale=.6](0,1.5)(.5,3)(2,1.5)(3.5,3)(4,1.5)
\end{pspicture}
}

\newcommand{\FIGLEMi}{
\begin{pspicture}[.4](0,-1)(4.5,3.5)
\psdots[dotscale=.6](0,0)(2,0)(1,1.5)
\psline[linewidth=.5pt,linestyle=dotted,dotsep=1pt](0,0)(1,.5)(2,0)
\psline[linewidth=.5pt,linestyle=dotted,dotsep=1pt](1,.5)(1,1.5)
\pscurve[linewidth=.5pt](1,1.5)(1.15,1.7)(.85,1.9)(1.15,2)(.85,2.2)(1.15,2.4)
\pscurve[linewidth=.5pt](0,0)(-.2,.1)(-.4,-.4)(-.6,0)(-.8,-.6)(-1,-.5)
\pscurve[linewidth=.5pt](2,0)(2.2,.1)(2.4,-.4)(2.6,0)(2.8,-.6)(3,-.5)
\put(.2,-.4){$\scriptstyle{2}$}
\put(1.6,-.4){$\scriptstyle{3}$}
\put(1.2,1.1){$\scriptstyle{1}$}
\put(-1.3,0){$\scriptstyle{B}$}
\put(2.8,.1){$\scriptstyle{C}$}
\put(.2,2){$\scriptstyle{A}$}
\end{pspicture}
}

\newcommand{\FIGLEMii}{
\begin{pspicture}[.4](0,-1)(4.5,3.5)
\psdots[dotscale=.6](0,0)(2,0)(1,1.5)
\psline[linewidth=.5pt,linestyle=dotted,dotsep=1pt](0,0)(1,.5)(2,0)
\psline[linewidth=.5pt,linestyle=dotted,dotsep=1pt](1,.5)(1,1.5)
\pscurve[linewidth=.5pt](1,1.5)(1.15,1.7)(.85,1.9)(1.15,2)(.85,2.2)(1.15,2.4)
\pscurve[linewidth=.5pt](0,0)(.05,-.25)(-.4,0)(-.4,-.5)(-.8,-.2)(-1,-.5)
\pscurve[linewidth=.5pt](2,0)(1.95,-.25)(2.4,0)(2.4,-.5)(2.8,-.2)(3,-.5)
\put(.2,-.4){$\scriptstyle{3}$}
\put(1.6,-.4){$\scriptstyle{2}$}
\put(1.2,1.1){$\scriptstyle{1}$}
\put(-1.3,0){$\scriptstyle{C}$}
\put(2.8,0){$\scriptstyle{B}$}
\put(.2,2){$\scriptstyle{A}$}
\end{pspicture}
}

\newcommand{\FIGPFi}{
\begin{pspicture}[.4](0,-1)(4.5,3.5)
\psdots[dotscale=.6](0,0)(2,0)(0,2)(2,2)
\psline[linewidth=.5pt,linestyle=dotted,dotsep=1pt](0,0)(0,2)
\psecurve[linewidth=.5pt,linestyle=dotted,dotsep=1pt](-1,1)(0,1)
         (1,1)(2,0)(2.5,-1)
\pscurve[linewidth=.5pt](0,0)(-.2,.1)(-.4,-.4)(-.6,0)(-.8,-.6)(-1,-.5)
\pscurve[linewidth=.5pt](0,2)(-.2,2.1)(-.4,2.4)(-.6,2)(-.8,2.6)(-1,2.5)
\pscurve[linewidth=.5pt](2,0)(2.3,-.1)(2.1,.4)(2.6,.6)(2.5,.8)
        (2.6,1)(2.5,1.2)(2.6,1.4)(2.1,1.6)(2.3,2.1)(2,2)
\put(1.4,2){$\scriptstyle{4}$}
\put(.2,-.4){$\scriptstyle{2}$}
\put(1.4,-.4){$\scriptstyle{3}$}
\put(-1.7,2){$\scriptstyle{A}$}
\put(-1.4,-.1){$\scriptstyle{B}$}
\put(2.8,.8){$\scriptstyle{C''}$}
\put(1.6,.8){$\scriptstyle{C'}$}
\put(.2,2){$\scriptstyle{1}$}
\end{pspicture}
}

\newcommand{\FIGPFii}{
\begin{pspicture}[.4](0,-1)(4.5,3.5)
\psdots[dotscale=.6](0,0)(2,0)(0,2)(2,2)
\psline[linewidth=.5pt,linestyle=dotted,dotsep=1pt](0,0)(0,2)
\psecurve[linewidth=.5pt,linestyle=dotted,dotsep=1pt](-1,1)(0,1)(1,1)
         (2,2)(2.5,3)
\pscurve[linewidth=.5pt](0,0)(-.2,.1)(-.4,-.4)(-.6,0)(-.8,-.6)(-1,-.5)
\pscurve[linewidth=.5pt](0,2)(-.2,2.1)(-.4,2.4)(-.6,2)(-.8,2.6)(-1,2.5)
\pscurve[linewidth=.5pt](2,0)(2.3,-.1)(2.1,.4)(2.6,.6)(2.5,.8)(2.6,1) 
        (2.5,1.2) (2.6,1.4)(2.1,1.6) (2.3,2.1)(2,2)
\put(1.4,2){$\scriptstyle{4}$}
\put(.2,-.4){$\scriptstyle{2}$}
\put(1.4,-.4){$\scriptstyle{3}$}
\put(-1.7,2){$\scriptstyle{A}$}
\put(-1.4,-.1){$\scriptstyle{B}$}
\put(2.8,.8){$\scriptstyle{C''}$}
\put(1.6,.8){$\scriptstyle{C'}$}
\put(.2,2){$\scriptstyle{1}$}
\end{pspicture}
}

\newcommand{\Gii}{
\begin{pspicture}[.4](0,-.2)(1,1)
\psline(0,0)(1,0)(.5,.85)(0,0)
\psdots(0,0)(1,0)(.5,.85)
\put(1.1,-.3){$\scriptstyle{2}$}
\put(-.2,-.3){$\scriptstyle{1}$}
\put(.45,1){$\scriptstyle{3}$}
\end{pspicture}
}
\newcommand{\Gia}{
\begin{pspicture}[.4](0,-.2)(1,1)
\psline(0,0)(1,0)(.5,.85)
\psdots(0,0)(1,0)(.5,.85)
\put(1.1,-.3){$\scriptstyle{2}$}
\put(-.2,-.3){$\scriptstyle{1}$}
\put(.45,1){$\scriptstyle{3}$}
\end{pspicture}
}
\newcommand{\Gib}{
\begin{pspicture}[.4](0,-.2)(1,1)
\psline(1,0)(.5,.85)(0,0)
\psdots(0,0)(1,0)(.5,.85)
\put(1.1,-.3){$\scriptstyle{2}$}
\put(-.2,-.3){$\scriptstyle{1}$}
\put(.45,1){$\scriptstyle{3}$}
\end{pspicture}
}
\newcommand{\Gic}{
\begin{pspicture}[.4](0,-.2)(1,1)
\psline(.5,.85)(0,0)(1,0)
\psdots(0,0)(1,0)(.5,.85)
\put(1.1,-.3){$\scriptstyle{2}$}
\put(-.2,-.3){$\scriptstyle{1}$}
\put(.45,1){$\scriptstyle{3}$}
\end{pspicture}
}

\def\lmn#1{\vadjust{\setbox1=\vtop{\hsize 12mm
  \parindent=0pt\baselineskip=9pt
  \rightskip=4mm plus 4mm#1}
  \hbox{\kern-12mm\smash{\raise .5ex\box1}}}}

\newtheorem{theorem}{Theorem}[section]

\newtheorem{lemma}[theorem]{Lemma}
\newtheorem{proposition}[theorem]{Proposition}
\newtheorem{definition}[theorem]{Definition}
\newtheorem{corollary}[theorem]{Corollary}
\newtheorem{remark}[theorem]{Remark}

\renewcommand{\det}{\operatorname{det}}
\newcommand{\Pf}{\operatorname{Pf}}
\newcommand{\adj}{\operatorname{adj}}

\def\Sy{\mathfrak S}
\def\P{\mathcal P}
\def\o{\mathfrak o}
\def\s{\sigma}
\def\LP{{\Lambda^{(p)}}}
\def\Z{{\mathbb Z}}
\def\R{{\mathbb R}}
\def\Rm{\mathcal{V}}
\def\X{\mathcal X}
\def\Y{{\mathrm Y}}

\begin{document}

\title[A New Matrix-Tree Theorem]{A New Matrix-Tree Theorem}


\author{Gregor Masbaum}
\address{Institut de Math{\'e}matiques de Jussieu (UMR 7586 CNRS),
Equipe `Topologie et g{\'e}om{\'e}trie alg{\'e}briques', Case 7012,
Universit{\'e} Paris VII, 75251 Paris Cedex 05,  France} 
\email{masbaum@math.jussieu.fr}

\author{Arkady Vaintrob}
\address{Department of Mathematics, University of Oregon, Eugene, OR
97405, USA} 
\email{vaintrob@math.uoregon.edu}

\begin{abstract} 
The classical Matrix-Tree Theorem allows one to list the spanning trees of a
graph by monomials in the expansion of the determinant of a certain matrix.   
We prove that in the case of three-graphs (that is,
hypergraphs whose edges have exactly \emph{three} vertices)
the spanning trees are generated by the Pfaffian of a suitably defined
matrix. This result can be interpreted topologically as an expression
for the lowest order term of the Alexander-Conway polynomial of an 
algebraically split link. We also prove some algebraic properties of our
Pfaffian-tree polynomial.  
\end{abstract}

\maketitle

\tableofcontents

\section{Introduction}

The classical Matrix-Tree Theorem of Kirchhoff
provides the following way of listing all the spanning 
trees in a graph.  
Consider a finite graph $G$ with vertex set $V$ and set of edges
$E$. Multiple edges between two vertices are allowed, and we denote by
$v(e) \subset V$ the set of endpoints of the edge $e$.  
If we label each edge $e\in E$  by a variable $x_e$, then a subgraph
of $G$ given as a collection of edges $S\subset E$ corresponds to the
monomial  $$ x_S=\prod_{e\in S}x_e.$$ 
Form a symmetric matrix $\Lambda=(\ell_{ij})$,  whose rows and columns 
are indexed by the vertices of the graph and entries given by
$$\ell_{ij}=-\sum_{{e\in E,} \atop {v(e)=\{i,j\} }} x_e, \text{\ \ if\
} i\ne j,   \text{\ \ \ and\ \ \ } 
 \ell_{ii}=\sum_{{e\in E,} \atop { i\in v(e)}}x_e .$$  
(This matrix arises in the theory of electrical networks and is called
sometimes the \emph{Kirchhoff matrix} of the graph.)
Since the entries in each row of $\Lambda$ add up to zero, the
determinant of this matrix vanishes and the determinant of the
submatrix $\Lambda^{(p)}$ obtained by deleting the $p$th row and
column of $\Lambda$ is independent of 
$p$.\footnote{In fact, all entries of $\adj(\Lambda)$, the matrix of 
cofactors of $\Lambda$, are equal. This can be seen as follows. 
By hypothesis, $\Lambda v=0$, where $v$ is the column vector with 
all entries equal to $1$. If $\Lambda$ has rank $\leq m-2$, where $m$ 
is the size of $\Lambda$, the matrix $\adj(\Lambda)$ is identically 
zero. Otherwise, $\Lambda$ has rank $m-1$, and its kernel is generated
by $v$. Now the formula $\Lambda \adj(\Lambda)=0$ shows that all 
columns of $\adj(\Lambda)$ are multiples of $v$. Since $\Lambda$ is symmetric,
so is $\adj(\Lambda)$, proving that all entries of $\adj(\Lambda)$ are equal.} 
This gives  a polynomial $\mathcal{D}_G$ in variables $x_e$ which is
called the Kirchhoff polynomial of $G$. 
The \emph{Matrix-Tree Theorem}~\cite[Theorem VI.29]{MTT} states that
non-vanishing  monomials appear in the polynomial $\mathcal{D}_G$  
with coefficient $1$ and correspond
to the spanning  trees
of $G$ ({\em i.e.} connected acyclic
subgraphs of $G$ with vertex set $V$). In other words,
\begin{equation}
  \label{eq:mtt}
  \mathcal{D}_G:=\det \Lambda^{(p)} = \sum_T x_T~,
\end{equation}
where the sum is taken over all the spanning trees in $G$.
\\[2pt]

For example, if $G$ is the complete graph $$K_3=\quad \Gii$$ with
vertices $1,2,3$  and edges $\{1,2\}, \{1,3\},$ and $\{2,3\}$,  we
have  
$$\Lambda = \left( \begin{array}{ccc}
x_{13}+x_{12} & -x_{12}& -x_{13}\\
-x_{12} &x_{12}+x_{23} & -x_{23}\\
-x_{13} &-x_{23} & x_{13}+x_{23}
\end{array}\right),
$$
and 
$$\mathcal{D}_G=\det\Lambda^{(1)}=\det \Lambda^{(2)}=\det
\Lambda^{(3)}=x_{12}x_{23}+x_{12}x_{13}+x_{13}x_{23},$$ 
which corresponds to the three spanning  trees of $K_3$ :\\[2pt]
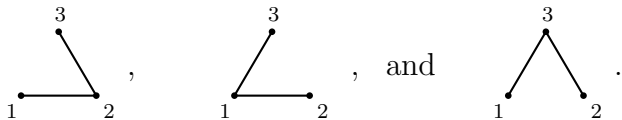
$$\Gia \quad , \qquad \quad \Gic\ \quad , \text{\ \ and} \  
\qquad \Gib \quad . $$  
\\[2pt]

\begin{figure}[h]
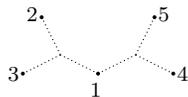

\begin{center}
\FIGEXi
\caption{\label{figexi} A spanning tree in the complete $3$-graph
$\Gamma_5$. It   has two edges, $\{1,2,3\}$ and $\{1,4,5\}$, 
and   contributes the term
  $ y_{123}\, y_{145}$ to $\P_5=P(\Gamma_5,\o_{can})$.}  
\end{center}
\end{figure}

In this paper we present 
an analog of this theorem for 
\emph{three-graphs} (or $3$-graphs). 
Edges of         a $3$-graph 
have three vertices and can be visualized
as triangles or  Y-shaped objects \MILN with the three vertices at their
endpoints.   We prove that the spanning trees of a $3$-graph 
can be generated by the terms in the expansion of the
 Pfaffian of a suitably defined \emph{skew-symmetric} matrix.  
(A  sub-$3$-graph of a $3$-graph $G$ is called a spanning tree  
if its vertex set coincides with that of $G$  and the ordinary graph  
obtained by gluing  together  Y-shaped objects   $ \MILN$ 
corresponding to the edges of $T$ is a tree. 
See  Figure~\ref{figexi} for an example.)
\medskip

Let us describe our result for the complete 
$3$-graph $\Gamma_m$ with the vertex set
$V(\Gamma_m)=\{1,2,\ldots,m\}$.  The edges of $\Gamma_m$  are the
three-element subsets $\{i,j,k\}$ of $V(\Gamma_m)$. 

Introduce variables $ y_{ijk}$, with $i,j,k \in V(\Gamma_m)$,
anti-symmetric in $i,j,k$, {\em i.e.}
\begin{equation}
  \label{eq:skew}
 y_{ijk}=-  y_{jik}=  y_{jki} \text{\ \ and \ \ }  y_{iij}=0~. 
\end{equation}
Consider the $m \times m$ matrix 
\begin{equation}
\label{Lambda}
\Lambda =(\lambda_{ij}), \ {1\leq i,j\leq m}, \quad \text{with} \quad
\lambda_{ij}=\sum_{k}  y_{ijk}. 
\end{equation}
 This matrix is skew-symmetric and its entries in each row 
add up to zero.  This implies that 
the determinant of $\LP$ is independent of $p$. (Here, as before, 
$\LP$ denotes the result of removing the $p$th row and column from
$\Lambda$.) 

For example, if $m=3$, 
we have
$$
\Lambda = \left( \begin{array}{ccc}
0 &  y_{123}&  y_{132}\\
 y_{213} & 0&  y_{231}\\
 y_{312} &  y_{321} & 0
\end{array}\right) 
= 
\left( \begin{array}{ccc}
0 &  y_{123}& - y_{123}\\
- y_{123} & 0&  y_{123}\\
 y_{123} & - y_{123} & 0
\end{array}\right)
$$
and $\det \Lambda^{(3)} = - y_{123} y_{213}= y_{123}^2$.

Since the submatrix $\LP$ is still skew-symmetric,
it has a \emph{Pfaffian}, $\Pf( \LP)$,  whose square is equal to $\det
\LP$.  (See Section \ref{pfaffians} for a 
review of Pfaffians and their properties.) 

It turns out that  $(-1)^{p-1}\Pf( \LP)$ does not depend on $p$
(see Lemma~\ref{4.1})  
which allows us to define the polynomial
\begin{equation}
 \label{fctP}
\P_m : = (-1)^{p-1} \Pf(\LP)~
\end{equation}
in variables $ y_{ijk}$.
We will call $\P_m$  the \emph{Pfaffian-tree polynomial} because of
its connections with spanning trees in $\Gamma_m$.

In the example with $m=3$ above,  one has   
$$\P_3=\Pf(\Lambda^{(3)})= y_{123}.$$ 

Note that since the determinant of 
a skew-symmetric matrix of odd size
is always zero, $\P_m=0$   if $m$ is even.

As in the case of ordinary graphs, the correspondence
\begin{equation*}
  \label{eq:sgr}
\text{variable \ }   y_{ijk} \quad  \mapsto \quad  \text{edge\ }
\{i,j,k\} \ \text{of \ } \Gamma_m
\end{equation*}
assigns to each monomial in $ y_{ijk}$  a 
sub-$3$-graph of $\Gamma_m$.
A remarkable property of the polynomial $\P_m$ is that the 
sub-$3$-graphs
of $\Gamma_m$ corresponding to its monomials are precisely
the spanning trees of $\Gamma_m$. In particular, if $m$ is odd,
the $3$-graph $\Gamma_m$ has no spanning trees (this, of course, can
be easily proved directly by a simple combinatorial argument). 

Note, however, that because of condition~(\ref{eq:skew}) the
correspondence between monomials and 
sub-$3$-graphs is not one-to-one. A sub-$3$-graph
determines a monomial only up to a sign. 

In order to express the signs of the monomials of $\P_m$ in terms of
spanning trees, we introduce a notion of orientation of $3$-graphs
(see Section~\ref{ori}).  
The $3$-graph $\Gamma_m$ has a canonical  orientation $\o_{can}$ 
given by the natural ordering of the vertices. 
If  $T$ is a spanning tree in $\Gamma_m$, this orientation allows us
to define unambiguously  a monomial $y(T,\o_{can})$ 
(see Section~\ref{enum}
for details) which is, up to sign, just the product of the variables
$ y_{ijk}$ over the edges of $T$.  
The sum of these monomials is the generating function for spanning
trees in $\Gamma_m$, denoted by $P(\Gamma_m,\o_{can})$.  
Our Pfaffian Matrix-Tree Theorem in the case of the complete $3$-graph
states, then, that this generating function is given by the Pfaffian-tree
polynomial $\P_m$:  
\begin{equation}\label{eq:pfmtt} 
P(\Gamma_m,\o_{can}) = \P_m
\end{equation}

For example, if $m=5$, we have
\begin{equation} \label{m=5}
\P_5=P(\Gamma_5,\o_{can}) =   y_{123}\, y_{145} \  \pm \ \ldots~,
\end{equation}
where the right-hand side is a sum of $15$  
terms corresponding
to the $15$ spanning trees of $\Gamma_5$.
In Section~\ref{enumex} we will explain how to write the terms of
$P(\Gamma_m,\o_{can})$ explicitly including signs. In the case $m=5$,
all spanning trees are obtained from each other by permutations and
the right-hand side of~(\ref{m=5}) can be written with signs as 
$$ {1 \over 8}\sum_{\s\in{\Sy}_5} (-1)^\s
 y_{\s(1)\s(2)\s(3)}
\, y_{\s(1)\s(4)\s(5)}~.$$

If we visualize the edges of $\Gamma_m$ as  Y-shaped
objects \MILN, then the spanning tree corresponding to the first term
of~(\ref{m=5})  will look like on Figure~\ref{figexi}.  

\medskip
This paper grew out of our work~\cite{MV2} on connections
between the Alexander-Conway polynomial of links in $\R^3$ and the 
Milnor higher linking numbers. 
In Section~\ref{sec:top} we give a brief overview of the relations
between the  matrix-tree theorems and invariants of links. 
In Section~\ref{sec:3graphs}, we formally introduce $3$-graphs, their
spanning trees, and orientations.  There we also deal with the issue
of signs and define the generating function $P(G,\o)$ for spanning
trees in an oriented $3$-graph $G$. 
In Section~\ref{pfaffians}, we review Pfaffians and prove several
properties we need. Our Matrix-Tree Theorem expressing $P(G,\o)$ as
the Pfaffian-tree polynomial $\P_G$ of $G$ is proved in
Section~\ref{nmt}. Finally in Section~\ref{sec:prop}, we establish
some interesting algebraic properties of the Pfaffian-tree
polynomial. In particular, we show that $\P_G$ satisfies a
contraction-deletion relation. We also prove a $3$-term relation  and
a $4$-term relation for the Pfaffian-tree polynomial $\P_m$ of the
complete $3$-graph $\Gamma_m$, and give a recursion formula for
$\P_m^2$ which can be used to identify this polynomial in other
contexts.  
\medskip

\noindent \emph{Acknowledgements.}
We would like to express our gratitude to institutions  
whose hospitality and support we enjoyed during our work on this
paper: 
the Max-Planck Institut f{\"u}r Mathematik in Bonn (A.V.), 
IH{\'E}S (A.V.), the University of Oregon (G.M.),  and the 
University of Paris VII 
(A.V.).  
We also thank  R.~Booth and A.~Borovik for valuable comments
and the journal referee for suggestions on improving the exposition.
Research of the second author was partially supported by NSF grant
DMS-0104397.

\section{Alexander-Conway polynomial, linking numbers, and Milnor
invariants} 
\label{sec:top} 

In this section we briefly discuss the topological motivation of our
work. It is not necessary for understanding the rest of the paper and
may be safely  skipped by a reader interested only in the
combinatorial aspect of our   results.
\medskip

Let $L$ be an oriented link in  $\mathbb{R}^3$. The best
known classical isotopy invariants of $L$ are the linking numbers
$\ell_{ij}(L)$ between the $i$th and $j$th components of $L$ and 
its Alexander-Conway polynomial 
(see {\em e.g.}~\cite{K}) 
$$\nabla_L = \sum_{i\geq 0} c_i(L) z^i \in \Z[z].$$

After the work of Hosokawa~\cite{Hw}, Hartley~\cite[(4.7)]{Ha},  
and Hoste~\cite{Ho}, it is known that the coefficients $c_i(L)$ of
$\nabla_L$ for an $m$-component link $L$ vanish when $i\leq m-2$ 
and the coefficient $c_{m-1}(L)$ depends only
on the linking numbers $\ell_{ij}(L)$ 
via the determinantal formula
\begin{equation}
  \label{eq:ho_det}
c_{m-1}(L)= \det\Lambda^{(p)},
\end{equation}
where $\Lambda=(\lambda_{ij})$ is the matrix formed by linking numbers 
$$\lambda_{ij}=
\begin{cases}
\ \ \ \ \ -\ell_{ij}(L), & \text{\ if \ $i\ne j$\ }  \\
\sum_{k\neq i} \ell_{ik}(L), & \text{\ if \ $i=j$}~.
\end{cases}
$$

Hartley and Hoste also gave a second expression for
$c_{m-1}(L)$ as a sum over trees: 
\begin{equation}\label{eq:ho_tree}
c_{m-1}(L)= \sum_{T} \prod_{\{i,j\}\in edges(T)} \ell_{ij}(L)~,
\end{equation}
 where $T$ runs through the maximal trees 
in the complete graph $K_m$ with vertices $\{1,2,\ldots,m\}$.

For example, if $m=2$ then $c_1(L)=\ell_{12}(L)$, and if $m=3$, then
$$c_2(L)=\ell_{12}(L)\ell_{23}(L) +
\ell_{12}(L)\ell_{13}(L)+\ell_{13}(L)\ell_{23}(L)~,$$ 
corresponding to the three spanning
trees of $K_3$. 

The equality of the expressions~(\ref{eq:ho_det})
and~(\ref{eq:ho_tree}) for  $c_{m-1}(L)$ follows from the classical
Matrix-Tree Theorem \eqref{eq:mtt} applied to the complete graph with
$m$ vertices. 

If the link $L$ is \emph{algebraically split}, {\em i.e.} all linking 
numbers  $\ell_{ij}$ vanish, then not only  $c_{m-1}(L)=0$, but, as
was proved by Traldi~\cite{Tr} and Levine~\cite{Le}, the next $m-2$
coefficients of $\nabla_L$ also vanish
$$c_{m-1}(L)=c_{m}(L)=\ldots=c_{2m-3}(L)=0.$$   

For an algebraically split oriented link $L$ with three components,
there exists an integer-valued isotopy invariant $\mu_{123}(L)$ called
the \emph{Milnor  triple linking number} 
                   (see~\cite{Mi}).
(This invariant is equal to $1$ 
for the standard Borromean rings.) For an algebraically split link $L$
with $m$ components, we thus have $m \choose 3$ triple linking numbers
$\mu_{ijk}(L)$ corresponding to the different $3$-component sublinks
of $L$.  Unlike the ordinary linking numbers, the triple linking
numbers  are anti-symmetric in their indices  
$$\mu_{ijk}(L)=-\mu_{jik}(L)=\mu_{jki}(L)~.$$

Levine~\cite{Le} (see also Traldi~\cite[Theorem~8.2]{Tr2})
 found an expression for the first non-vanishing coefficient
$c_{2m-2}(L)$ of $\nabla_L$ for an algebraically split link in terms
    of  Milnor triple linking numbers
\begin{equation}
  \label{eq:lev}
c_{2m-2}(L)=\det\LP~,
\end{equation}
where $\Lambda=(\lambda_{ij})$ is the $m\times m$ skew-symmetric
matrix with entries given by 
$$
\lambda_{ij}=\sum_{k} \mu_{ijk}(L)
$$ 
 (cf.~(\ref{Lambda})).

This formula is analogous to the 
determinantal expression~(\ref{eq:ho_det}).
One of the starting points of the present paper was an attempt to find  
an analogue of the sum over trees formula~(\ref{eq:ho_tree}) for
algebraically split links. 
As a corollary of the determinantal formula~(\ref{eq:lev})
and our Matrix-Tree Theorem for complete $3$-graphs we obtain a
combinatorial expression for $c_{2m-2}$  as  
(the square of) a sum over trees. 
 
\begin{theorem}  \label{thm:alex-trees}
If $L$ is an algebraically split link, then
  \begin{equation}
    \label{eq:sqtrees}
c_{2m-2}(L)= P_m(L)^2~, 
  \end{equation}
where $P_m(L)$ is the spanning tree generating function
  $P(\Gamma_m,\o_{can})$ evaluated at $y_{ijk}=\mu_{ijk}(L)$. 
\end{theorem}
\proof This follows from (\ref{eq:lev}) and (\ref{eq:pfmtt}), since
$\det  \LP= \Pf( \LP)^2$. \eeproof 

\medskip

Using the theory of finite type invariants, we give in~\cite{MV2}
 another, direct, proof of formula (\ref{eq:sqtrees}) for the
 coefficient $c_{2m-2}(L)$ which does not use the determinantal 
formula~(\ref{eq:lev}). In fact, the proof  in~\cite{MV2} together  
 with our Matrix-Tree Theorem for   $3$-graphs provides an alternative 
 proof of formula~(\ref{eq:lev}).  

The following table summarizes the analogy which was our guiding
principle in this work: 

 \medskip
\begin{center}
 {\small
\begin{tabular}[c]{|c|c|}
\hline
Linking numbers & Milnor's triple linking numbers\\ 
 Edges of ordinary graphs& (Oriented) Edges of $3$-graphs \\
Formulas~(\ref{eq:ho_det}) and~(\ref{eq:ho_tree})& 
Formulas~(\ref{eq:lev}) and~(\ref{eq:sqtrees})\\
The classical Matrix-Tree Theorem& The Pfaffian Matrix-Tree Theorem \\
\hline
 \end{tabular}}
\end{center}

\section{Three-graphs}\label{sec:3graphs}

\subsection{Basic definitions} \

\emph{Three-graphs}
(or, more officially, \emph{$3$-uniform
  hypergraphs}, see~\cite{Ber}) are 
analogues of  graphs.  
The only difference is that the edges of a $3$-graph are `triangular',
 {\em i.e.\/} they  have three vertices,
while edges of ordinary graphs have only two vertices.
For our purposes, it will not be necessary to consider degenerate
 edges, {\em i.e.\/}  edges with less than $3$ vertices (they are
analogous to loop   edges for ordinary graphs), and so we will use the
following definition. 

\begin{definition}\rm  
A \emph{three-graph}  (or $3$-graph)
is a triple $G=(V,E,v)$ where $V=V(G)$ and
  $E=E(G)$ are   finite  sets, 
whose elements will be called, respectively, \emph{vertices} and
  \emph{edges} of $G$, and $v$ is a map from $E$ to the set 
of three-element subsets of $V$.
For an edge $e\in E$ the elements of $v(e)$ are called the \emph{vertices}
or the \emph{endpoints} of $e$.  
\medskip

A  $3$-graph $G'=(V',E',v')$ is called 
a  {\em sub-$3$-graph\/}  of  $G=(V,E,v)$
if $V'\subset V$, $E'\subset E$ and $v'=v_{|E'}$. 
If $V'=v(E')$, we say that $G'$ is the 
sub-$3$-graph  of $G$ generated by the subset of edges $E'$.
\end{definition}

Note that by replacing ``three-element'' with ``two-element''
in the above definition, we recover exactly the 
definition of a (finite) graph without loops, but, possibly, with
multiple edges.  

Visually, an edge of a $3$-graph can be thought of as a triangle or 
as a Y-shaped object \MILN with the three vertices at its endpoints.  

In this paper we adopt the latter point of view and, accordingly,
we define the \emph{topological realization} of a $3$-graph by
gluing together all the Y's corresponding to its edges.

More formally,
the topological realization $|G|$ of a $3$-graph $G$ is a one-dimensional
cell complex obtained by taking one $0$-cell for every element 
of $V(G) \cup E(G)$ and then gluing  a $1$-cell  for every pair
$(v,e)\in V(G) \times E(G)$ such that 
$v\in v(e)$.  

Note that $|G|$ is the same as the topological realization of the 
\emph{bipartite graph} naturally associated 
with
the $3$-graph  $G$.

\subsection{Trees} \

\begin{definition}\rm
A $3$-graph $G$ is called a \emph{tree} if its topological realization
$|G|$ is connected and simply connected.  
({\em I.e.} the bipartite graph formed by the 
Y-shaped objects \MILN corresponding to 
the edges of $G$ is a tree.) 

\end{definition}

The following proposition states
that this definition is equivalent to the standard 
definition of trees for hypergraphs (see~\cite{Ber})
and that a $3$-graph with an even number of vertices cannot be a tree. 

A  \emph{path\/} in a $3$-graph $G$
is a sequence $v_0,e_1,v_1,\ldots,v_{n-1},e_n,v_n$ of vertices $v_i$ 
and edges $e_i$ of $G$ such that $v_i \in v(e_{i+1})$ for $i=0,\ldots,n-1$,
and $v_i \in v(e_i)$ for $i=1,\ldots,n$. A path is called a \emph{cycle\/}
if $v_0=v_n$ and $e_i \ne e_j$ for $i\ne j$. A $3$-graph $G$ is called
\emph{connected\/} if for any two vertices $v$ and $u$ of $G$ there exists a
path in $G$ which begins in $v$ and ends in $u$. (In particular, a $3$-graph
with only one vertex and no edges is connected.)

\begin{proposition}\label{Eul}\

(i) A $3$-graph is a tree if and only if it is connected and has no cycles.

(ii) If a $3$-graph with $m$ vertices is a tree, then $m$ is odd and
the number of edges is equal to $(m-1)/2$.
\end{proposition}

Both statements can be easily proved by induction on the number of edges.

A  sub-$3$-graph  $T$ of a $3$-graph $G$ is called a \emph{spanning tree} if
$T$ is a tree  and  $V(T)=V(G)$.
By the above proposition, only $3$-graphs with odd
number of vertices may have spanning trees.
\medskip

In the next section we will need the following
characterization of trees.

\begin{proposition} \label{prop:trees}
Let $T$ be a $3$-graph with $n$ edges and $2n+1$ vertices.
Fix an ordering of the vertex set $E(T)$, and for each edge $e\in E(T)$,
choose a  (non-trivial) cyclic permutation $\sigma(e)$ of the
three-element  set $v(e)$.
Then $T$ is a tree if and only if the product 
\begin{equation}  \label{eq:prodcycl}
\sigma(T)=\prod_{e\in E(T)}\sigma(e)
\end{equation}
is a cyclic permutation of the vertex set $V(T)$. 
In particular, if the product in some order is a cyclic 
permutation of $V(T)$,  then the same is true for any 
other order as well. 
  
\end{proposition}
\proof
First, assume that $T$ is a tree.
We will prove that $\sigma(T)$ is a cycle by induction on the number
of edges in $T$. In the case $n=1$, 
the statement is a tautology. If $n\geq 2$, number the edges
$e_1,e_2,\ldots, e_n\in E(T)$ according to the order they appear in
the product~\eqref{eq:prodcycl}, so that  
$$
\sigma(T)=\sigma(e_1)\sigma(e_2)\ldots\sigma(e_n)~.
$$  
Since $T$ is a tree, it has an edge $e_k$ which has only one common
vertex with the sub-$3$-graph  $T'$ generated by the remaining edges 
(in which case, $T'$ is itself a tree). Since a conjugate of
any cyclic permutation  is again a cyclic permutation, the equation 
\begin{equation} \label{eq:product}
\sigma(e_1)\sigma(e_2)\ldots\sigma(e_n) =
\pi^{-1}(\sigma(e_k)\ldots\sigma(e_n) 
\sigma(e_1) \ldots \sigma(e_{k-1})) \pi~,
\end{equation}
where $\pi=\sigma(e_k)\ldots\sigma(e_n)$,
shows that we can assume that $k=1$, {\em i.e.}\ that the 
sub-$3$-graph $T'$ of $T$ generated by the edges $e_2,e_3,\ldots,e_n$ is 
also a tree. Therefore, by induction, the permutation
$$
        \sigma(      {T}')=\sigma(e_2)\sigma(e_3)\ldots\sigma(e_n)
$$
 is a cyclic permutation of the set $V(T')$. This implies that
the permutation
        $\sigma({T})=\sigma(e_1)\sigma({T}')$ 
is a cycle, since it is a product of two cyclic permutations having
only one common element.

Conversely, assume that $T$ is not a tree. Then it is disconnected, because a 
connected $3$-graph with $n$ edges can have $2n+1$ vertices only if it is a
tree. 
Therefore, the permutation $\sigma(T)$  cannot be a $2n+1$-cycle, 
since it is equal to the product of several commuting permutations
corresponding to the different components of $T$.
\eeproof

\begin{remark}
\rm
A similar description of trees exists for ordinary graphs. 
Namely, every edge of a graph $G$ determines a transposition (a
two-cycle) of the vertex set $V(G)$; and a graph with $m$ vertices and 
$m-1$ edges is a tree  if and only if the product (taken in any order) 
of the corresponding $m-1$ transpositions is an $m$-cycle.  
\end{remark}

\subsection{Orientations}\label{ori} \

In order to keep track of signs, we 
introduce a notion of \emph{orientation} as
follows.    
\begin{definition}\rm
An \emph{orientation} of a finite set $\X$ is an ordering
 of $\X$ up to even permutations. An \emph{orientation} of a   
$3$-graph $G$ is an orientation of its vertex set $V(G)$. 
An \emph{orientation} of an edge $e\in E$ is an orientation of its
  vertex set  $v(e)$.  
\end{definition}

Note that an orientation of an edge $e\in E(G)$ is the same as an
orientation of the sub-$3$-graph of $G$ generated by that edge.

\begin{remark} \rm
The term `orientation' here is justified because an orientation of 
a set $\X$ induces an orientation of the vector space
$\mathbb{R}^{\X}$.
\end{remark}

If a $3$-graph $G$ has an odd number, $m$, of vertices,
then an orientation of $G$ 
can also be specified by a  {\em cyclic ordering}  of 
$V(G)$, {\em i.e.} an ordering up to cyclic permutation.
(This is because an $m$-cycle is an even permutation if $m$ is odd.)
We will usually  write such an ordering as a cyclic permutation.
In particular, an orientation of an edge of a $3$-graph
with vertex set $\{i,j,k\}$ 
will be indicated by choosing one of the
two three-cycles $(ijk)$ or  $(jik)$.   

Note that the orientations given by two $m$-cycles
$\sigma$ and $\sigma'$ are the same if and only if 
$\sigma'=s \sigma s^{-1}$ where
$s$ is an even permutation of $V(G)$. 
\medskip

The following is a key construction needed for the definition of the
generating function of spanning trees in a $3$-graph. 
Let $\tilde T$ denote the data consisting
of a tree $T$ together with a choice of orientation for every edge of 
$T$. For each edge $e\in E(T)$, denote by $\sigma(e)$ the cyclic permutation
of the three-element set $v(e)$ induced by the orientation of this
edge in $\tilde{T}$. 
If we choose an ordering of the set $E(T)$, then by
Proposition~\ref{prop:trees} 
the product 
\begin{equation}
  \label{eq:orprod}
\sigma(\tilde{T})=\prod_{e\in E(T)}\sigma(e)
\end{equation}
is a cyclic permutation of the set $V(T)$ and, therefore,
gives a cyclic ordering of $V(T)$.
Since $T$ has an odd number of vertices, this 
cyclic ordering induces an orientation of $T$.
As we will show now, this orientation does not depend on the 
choice of ordering of the edges and we denote it by $\o(\tilde T)$.

\begin{proposition} \label{signs} 
The orientation $\o(\tilde T)$ is well-defined. Moreover, 
$\o(\tilde T)$ changes sign  
whenever the orientation of an edge in $\tilde T$ is reversed.   
\end{proposition}

\proof
The proof that the orientation $\o(\tilde T)$ is well-defined,
\emph{i.e.} it does not depend on the order of factors
in~(\ref{eq:orprod}), is similar to the proof of 
Proposition~\ref{prop:trees}.
Indeed, if we change the order of the factors cyclically 
then $\sigma(\tilde{T})$ is replaced by 
$\pi^{-1}\sigma(\tilde{T})\pi$  (see~(\ref{eq:product})),
where $\pi$ is an even permutation (a product of several 
three-cycles $\sigma(e_i)$), and thus, the orientation
will not change. This implies that, as in the proof
of Proposition~\ref{prop:trees}, we may assume that the 
first factor $\sigma(e_1)$ in~(\ref{eq:orprod}) comes from an edge
that  has only one common vertex with the 
sub-$3$-graph $T'$ generated by
the remaining  edges. But in this case, $T'$ is a tree and by
induction we see that the orientation $\o(\tilde{T}')$, and therefore
$\o(\tilde T)$ too, is well-defined. 

Let us now prove that changing the orientation of one edge $e$ in
$\tilde T$ reverses the orientation $\o(\tilde T)$.

Denote the vertices of the tree $T$ by $V(T)=\{1,2,\ldots,m\}$ 
so that $v(e)=\{1,2,3\}$. Let 
$T_i$ be the maximal subtree of $T$ which contains the vertex $i$ 
but not the edge $e$.\footnote{It may happen that $T_i$ consists of 
just the vertex $i$ and no edge. Note that with our definition, 
such a $3$-graph is a tree.}
The data $\tilde{T}$ induces orientations of each of the
subtrees $T_1$, $T_2$ and $T_3$,
which we represent, respectively, by cyclic permutations
$(1A)$, $(2B)$ and $(3C)$, 
where $A$, $B$ and $C $ are ordered disjoint sets (some of which may
be empty)  with (unordered) union $\{4,5,\ldots ,m\}$.

\begin{figure}[h]
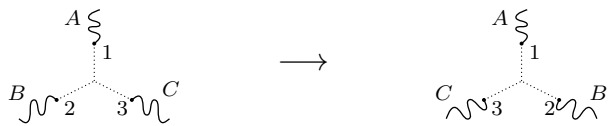

\begin{center}
 \FIGLEMi \ \ \ \ 
$\longrightarrow$ 
\ \ \ \ \ \ \ \ \ \ \  \ \
 \FIGLEMii 
\caption{\label{figlem} An admissible embedding of $\tilde T$ with the
orientation of the edge  $\{1,2,3\}$ reversed can be obtained by
exchanging the subtrees at the vertices $2$ and $3$. 
} 
\end{center}
\end{figure}

If the orientation of the edge $e$ in $\tilde{T}$ is
$(123)$, then the orientation $\o(\tilde T)$ is represented by
the product of cycles\footnote{
Our convention is that $\sigma_1\sigma_2$ means first apply $\sigma_2$
then $\sigma_1$.}
\begin{equation*}
\sigma=(1A)(123)(2B)(3C)=(A12B3C)~.
\end{equation*}
If we change the orientation of $e$ to $(132)$, then
the new orientation of $T$ is represented by
\begin{equation*}
\sigma'=(1A)(132)(2B)(3C)=(A13C2B)
\end{equation*}
(see Figure~\ref{figlem} for an illustration).  
Therefore,
$$\sigma'=s \sigma s^{-1}~,$$
where $s$ is an odd permutation,
since both $B$ and $C$ have 
an even number of elements.
Thus, $\sigma$ and $\sigma'$ define opposite orientations 
of the tree $T$. 

This completes the proof of the proposition. \eeproof
\medskip

In concrete situations it is often 
convenient to specify  orientations of the edges of trees by
embedding them 
into the plane. This leads to the following
pictorial way of computing the orientation $\o(\tilde{T})$.

If we embed the tree $T$ (or, more precisely, its topological realization
$|T|$)  into the plane, 
then
every  edge   of $T$ acquires
an orientation induced from the standard (counterclockwise) 
orientation of the plane.
We call such an embedding 
\emph{admissible} with respect to the data $\tilde T$ 
if for every edge
its induced 
orientation coincides with the one specified by $\tilde T$. 
By induction on the number of edges, it
is easy to see that 
for every $\tilde T$ there exists an admissible embedding.

\begin{figure}[h]
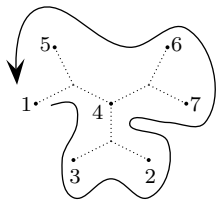

\begin{center}
\FIGEX
\caption{\label{figex} An admissible embedding of $\tilde T$ having
three oriented edges $(145)$, $(243)$, and $(476)$. 
Orderings representing $\o(\tilde T)$ include, for 
  example, $1432765$ and  $1327645$. They can be read off going 
  counterclockwise around the diagram, as indicated by the solid
arrow.}  
\end{center}
\end{figure}

Given an admissible embedding of $T$,  we can obtain an orientation
of $V(T)$ as follows. Go around the embedded tree  in the 
counterclockwise order
and write down the \emph{cyclic\/} sequence
of vertices in the order they were visited. 
Some vertices will be visited more than once, but when
we remove extra vertices from this sequence (in an arbitrary way) we
obtain a cyclic ordering of the set $V(T)$ (see an example on
Figure~\ref{figex}).  It is easy to see that the orientation
of $T$ given by this ordering is independent of the choices 
involved and coincides with the orientation $\o(\tilde T)$.   

As an example of this approach, Figure~\ref{figlem}
gives a pictorial illustration of the above proof of the sign change
property of the orientation $\o(\tilde{T})$.  

\subsection{The tree generating function} \
\label{enum}

We now introduce a function 
which generates spanning trees in a $3$-graph. 

Given a $3$-graph $G=(V,E,v)$, we denote by $\tilde E=\tilde E(G)$ the
set of oriented edges of $G$. Thus, an element $\tilde e$ of $\tilde
E(G)$ is a pair consisting of an edge $e\in E(G)$ and an orientation
of $e$. We denote by $-\tilde e$ the edge with the opposite
orientation.  The assignment $\tilde e \mapsto -\tilde e$ is a
fixed-point 
free involution on $\tilde E(G)$, with quotient $E(G)$.

To every oriented edge $\tilde e \in \tilde E(G)$, we associate an 
indeterminate $ y_{\tilde e}$ with the relation
$$ y_{-\tilde   e}=- y_{\tilde e}.$$

\begin{definition} {\em Let $T$ be a $3$-graph which is a tree, and
let $\o$ be an orientation of $T$. We define the monomial $y(T,\o)$ by
the formula  
 $$y(T,\o)={\o(\tilde T)\over \o}\prod_{e\in E(T)}  y_{\tilde e}~,$$
where $\tilde T$ is any choice of orientations for the edges (giving a
lift $\tilde e\in \tilde E(G)$ of every edge $e\in E(G)$), $\o(\tilde
T)$ is the orientation defined in Proposition~\ref{signs}, and
${\o(\tilde T)/ \o}$ is the sign relating the two orientations.}  
\end{definition}

If one changes the orientation of an edge $e$ in  $\tilde T$, then by
Proposition~\ref{signs} the orientation  $\o(\tilde T)$ picks up a
minus sign which is cancelled by the change of sign for $ y_{\tilde
e}$. This shows that  $y(T,\o)$ is well defined. 

\begin{definition} {\em Let $G$ be a $3$-graph, and let $\o$ be an
orientation of $G$. The} generating function for spanning trees {\em
  of $G$  is  
  \begin{equation}
    \label{eq:2x}
  P(G,\o)=\sum _T y(T,\o)~,  
  \end{equation}
where $T$ runs through the spanning trees of $G$.}
\end{definition}

Note that $P(G,-\o)=-P(G,\o)$. Thus, only the sign of the generating
function depends on the orientation.

\subsection{The tree generating function of the complete $3$-graph}
\label{enumex}
 \ 

Let us now give a more explicit combinatorial expression for the tree
generating function \eqref{eq:2x} of a   complete $3$-graph.  

By definition, the {\em complete $3$-graph}
$\Gamma_m$ has vertex set
$$V(\Gamma_m)=\{1,2,\ldots,m\},$$ and has exactly one  edge $e$ for
every unordered triple $\{i,j,k\}$ of vertices. Every 
cyclic permutation $(ijk)$ determines an \emph{oriented} edge $\tilde
e$. Therefore we can identify  the indeterminates 
$ y_{\tilde e}$ with  indeterminates $ y_{ijk}$ which are totally
antisymmetric 
in their indices (as in~(\ref{eq:skew})).  

We denote by $\o_{can}$ the orientation of $\Gamma_m$ given by the
natural ordering of $V(\Gamma_m)$. If $m$ is even, then $\Gamma_m$ has
no spanning trees. Therefore let us assume $m$ is odd. 
Put $d=(m-1)/2$.  

If $\s_1,\ldots, \s_{d}$ are $3$-cycles in $\Sy_m$, 
we set
$$
\sigma=\prod_i \sigma_i\in \Sy_m
$$
and define
$$
\varepsilon(\s_1,\ldots, \s_{d}) 
=
\begin{cases}
0 & \text{if $\sigma$ is not an $m$-cycle} \\
(-1)^s & \text{if $\sigma=(s(1)\ldots s(m))$ is an $m$-cycle,} \\
\end{cases}
$$
where $(-1)^s$ is the sign of the permutation $s\in \Sy_m$.
Notice that in the case when $\sigma$ is an $m$-cycle,
$s$ is defined by the condition
$$
\sigma = (s(1)\,s(2)\ldots s(m)) = s (1\,2 \ldots m) s^{-1}
$$
only up to powers of the standard cycle $(1\,2\ldots m)$.
However, since $m$ is odd, the sign $(-1)^s$ is well-defined.

If $T$ is a spanning tree on $\Gamma_m$, the associated monomial is 
$$y(T,\o_{can}) = \varepsilon\bigl( (i_1j_1k_1),\ldots,
(i_dj_dk_d)\bigr) \prod_{\alpha=1}^d y_{i_\alpha j_\alpha k_\alpha}$$ 
where  $T$ is given by a collection of (unoriented) edges
$e_\alpha=\{i_\alpha, j_\alpha, k_\alpha\}$, $1\leq \alpha \leq d$.  

Since by Proposition~\ref{prop:trees} the sign
$\varepsilon\bigl( (i_1j_1k_1),\ldots, (i_dj_dk_d)\bigr)$ is
zero if $T$ is not a spanning tree, we can write  
\begin{equation*}
  P(\Gamma_m,\o_{can}) =\sum \varepsilon\bigl( (i_1j_1k_1),\ldots,
  (i_dj_dk_d)\bigr) \prod_{\alpha=1}^d y_{i_\alpha j_\alpha k_\alpha}  
\end{equation*}
where the sum is taken over all monomials   $\prod_\alpha y_{i_\alpha
j_\alpha k_\alpha}$ of degree $d$,  with the convention that monomials
which differ only by changing the order of indices in some of the
variables  are taken only once. Alternatively, 
since $\varepsilon \ne 0$ implies that the triples
$(i_\alpha , j_\alpha , k_\alpha), \alpha=1,\ldots,d$
are distinct, we can write  
\begin{equation*}
  P(\Gamma_m,\o_{can}) =\frac{1}{6^d d!} \sum \varepsilon\bigl
  ( (i_1j_1k_1),\ldots, (i_dj_dk_d)\bigr) \prod_{\alpha=1}^d
  y_{i_\alpha j_\alpha k_\alpha}  
\end{equation*} where the sum is now over all $3d$-tuples of indices
  $i_1,j_1, \ldots k_{d}$. 
\medskip

Some interesting algebraic properties of the tree generating function
$P(\Gamma_m,\o_{can})$ will be given in Section \ref{sec:prop}.  

\section{Pfaffians}\label{pfaffians}

In this section, we review the definition and the main properties of 
Pfaffians (see {\em e.g.}~\cite[Chapter XIV.\S 10]{La} and
\cite[Chapter 9.5]{Pf} for more details). 
Let 
$$A=(a_{ij}), \ {1\leq i,j\leq 2n}, \ \ \ \ a_{ij}=-a_{ij}$$  
be 
a skew-symmetric
matrix.
One way to define its Pfaffian, $\Pf (A)$, is 
as follows. 
Associate to $A$ the $2$-form   
$$\alpha = \sum_{i<j} a_{ij}\  e_i \wedge e_j = {1\over 2} \sum_{i,j}
a_{ij}\ e_i \wedge e_j
$$ 
where $e_1, \ldots, e_{2n}$ is a basis of  
one-forms on a vector space of dimension $2n$.
Then  
\begin{equation} {\alpha^n \over n!} = \Pf(A)\  e_1\wedge e_2\wedge
\ldots \wedge e_{2n}\label{pf1} 
\end{equation}

For example, if $n=1$, then $\Pf(A)=a_{12}$, and if $n=2$, then
$\Pf(A)=a_{12} a_{34} - a_{13} a_{24} +a_{14} a_{23}$.  

Equation~(\ref{pf1}) implies the following 
  explicit
  formula
for $\Pf(A)$. 
\begin{eqnarray} 
\Pf(A) &= & {1\over 2^n n!} \sum_{\s\in \Sy_{2n}}\, (-1)^\s
  \,a_{\s(1)\s(2)}\, a_{\s(3)\s(4)} \cdots a_{\s(2n-1)\s(2n)} 
\nonumber \\
&=& \sum_{\s\in \Sy_{2n}\atop{
\s(1)<\s(3)<\ldots <\s(2n-1)\atop
\s(2i-1)<\s(2i)\  \text{for all $i$}}} (-1)^\s \,a_{\s(1)\s(2)}\,
a_{\s(3)\s(4)} \cdots a_{\s(2n-1)\s(2n)} \label{pf3} 
\end{eqnarray}

\medskip

We will need the following \emph{row development formula}. Let
  $A^{(i,j)}$   denote the $(2n-2) \times (2n-2)$-matrix 
  obtained by removing the $i$th and $j$th row and column from $A$. Then    
\begin{equation} \label{pf6} 
\Pf(A)=a_{12}\,\Pf(A^{(1,2)}) -a_{13}\,\Pf(A^{(1,3)}) +\ldots
  +a_{1,2n}\,\Pf(A^{(1,2n)})~.   
\end{equation} 
This formula can be deduced from (\ref{pf3}) by noticing
that $\s(1)=1$ for every permutation $\sigma$ appearing in the sum
in~(\ref{pf3}).   
\medskip

Here are two more standard properties of the Pfaffian:
 \begin{equation} 
  \Pf(A)^2=\det A \label{pf4}~,
\end{equation}  
\begin{equation} \label{pf5}
\Pf(S^T AS)=\det (S)\,\Pf(A)~.
\end{equation}

We omit the proofs, but notice that if $\det S\neq 0$, then 
$ S^TAS $ 
is just the matrix of the $2$-form $\alpha$ in a 
new basis, which allows one to deduce~(\ref{pf5}) from~(\ref{pf1}).  
\medskip

Finally, we need the following lemma.

\begin{lemma}\label{4.1} 
Let $A=(a_{ij})_{1\leq i,j\leq 2n+1}$ be   
 a skew-symmetric
$(2n+1)\times(2n+1)$ matrix  such that 
$$\sum_i a_{ij}=0$$ 
for all $j$. 
Then $$(-1)^{p-1}\, \Pf(A^{(p)})$$ is independent of $p$.  
\end{lemma}

\proof  
Define a bilinear form $\alpha$ on a $2n$-dimensional vector space
       with basis $b_1,b_2,\ldots, b_{2n}$ by
$\alpha(b_i,b_j)=a_{ij}$. 
Thus, $\alpha$ is the two-form 
associated to the matrix $A^{(2n+1)}$. 
Set 
$$b'=-\sum_{i=1}^{2n} b_i$$ 
and observe  that
$$\alpha(b',b_j)=-\sum_{i=1}^{2n} a_{ij}=a_{2n+1,j}~.$$ 
Thus, the matrix of the form $\alpha$ in the basis $b_1,b_2,\ldots,
b_{2n-1},b'$ is precisely $A^{(2n)}$. 
Since the corresponding base-change transformation
has determinant $-1$, it follows from 
(\ref{pf5}) that $$ \Pf(A^{(2n)})= - \Pf(A^{(2n+1)})~.$$ 
More generally, $A^{(p)}$ is the matrix of $\alpha$ in the basis
$$b_1,b_2,\ldots, b_{p-1}, b_{p+1}, \ldots,  b_{2n-1},b'~,$$ 
which is related to the standard basis by a base change of determinant
$(-1)^{p-1}$,  
which proves the lemma. 
\eeproof

\section{The Pfaffian Matrix-Tree Theorem}\label{nmt}

We now state and prove a formula expressing the generating function
for spanning trees in terms of a Pfaffian. 
Let $G$ be a $3$-graph with vertices 
numbered
from $1$ to $m$. 
Define the matrix $$\Lambda(G) =(\lambda_{ij})_{1\leq i,j\leq m}$$  as 
follows. The diagonal terms $\lambda_{ii}$ are zero. If $i\neq j$,
then    
\begin{equation*} 
\lambda_{ij}=\sum_{\tilde e}  y_{\tilde e}~,
\end{equation*}
where $\tilde e$ runs through the oriented edges $\tilde e\in \tilde
E(G)$ such that $i$ and $j$ are vertices of $\tilde e$, and the
orientation of $\tilde e$ 
is represented by the cyclic ordering $(ijk)$, where $k$ denotes the
third vertex of $\tilde e$. Since $ y_{-\tilde e}=- y_{\tilde e}$, the
matrix $\Lambda(G)$ is   
skew-symmetric.

\begin{definition}\rm
Let $G$ be a $3$-graph 
with vertices numbered from $1$ to $m$.  
The \emph{Pfaffian-tree polynomial} of $G$ is 
\begin{equation}
  \label{eq:pf-tree-pol}
\P_G =  (-1)^{p-1} \Pf(\Lambda(G)^{(p)})~,
\end{equation}
where $\Lambda(G)^{(p)}$ 
       for $p=1,2,\ldots,m$ 
is the matrix obtained by removing the $p$th row
and column from  $\Lambda(G)$. 
\end{definition}

Note that the the right-hand side of~(\ref{eq:pf-tree-pol}) is
independent of $p$ by Lemma~\ref{4.1}. 
\bigskip

\noindent{\bf Example.} In the case of the complete $3$-graph
$\Gamma_m$, we can write the variables $y_{\tilde e}$ of
$\P_{\Gamma_m}$ as $ y_{ijk}$, where  
the $ y_{ijk}$'s are 
totally antisymmetric 
with respect to their indices. Then the matrix $\Lambda(\Gamma_m)$ is
precisely  the matrix $\Lambda$
in~(\ref{Lambda}), so that the Pfaffian-tree polynomial
$\P_{\Gamma_m}$ is equal to the polynomial $\P_m$ defined in
\eqref{fctP}. 

\medskip

The name `Pfaffian-tree' is justified by the following theorem, which
is the main result of this paper.  

\begin{theorem} \label{mainth}
Let $G$ be a $3$-graph with vertices
numbered from $1$ to $m$.
 Then the generating function of  
spanning trees in  $G$ is equal to the Pfaffian-tree polynomial of $G$ 
\begin{equation} \label{eq:main}
P(G,\o_{can})= \P_G~.
\end{equation} 

Here, $\o_{can}$ is the canonical orientation
of $G$ determined by the 
ordering of the vertices. 
\end{theorem}

\proof
We will assume that $\P_G$ is given by formula~(\ref{eq:pf-tree-pol}) 
with $p=1$. We may assume $m$ is odd, since otherwise both sides of
(\ref{eq:main}) are zero.   

First, we will show that it is enough to prove the 
theorem in the case of the
complete $3$-graph $\Gamma_m$. Indeed, let us write the indeterminates
for $\Gamma_m$ as $ y_{ijk}$, where 
the $ y_{ijk}$ are 
totally antisymmetric 
in their 
indices. Then the generating function for $G$ is obtained from the one
for $\Gamma_m$ by the substitution 
\begin{equation*} \label{eq:1}
 y_{ijk}=\sum_{\tilde e}  y_{\tilde e}~,
\end{equation*}
where $\tilde e$ runs through the oriented edges 
$\tilde{e}=(e,\o)$ 
of $G$  such that $i,j$ and $k$ are the vertices of 
$e$, and the orientation $\o$ of $\tilde e$  
is represented by the cyclic permutation $(ijk)$. The same
substitution applied to the matrix $\Lambda=\Lambda(\Gamma_m)$ yields
the matrix $\Lambda(G)$. Thus, if the theorem holds 
for $\Gamma_m$, then it holds for any $G$. 

Let us now prove the result for $G=\Gamma_m$. We denote the spanning
tree generating function $P(\Gamma_m,\o_{can})$ by $P_m$.  Note that
both $P_m$ and the Pfaffian-tree polynomial $\P_m$ are polynomials in
the indeterminates $ y_{ijk}$.  

 We will prove that
\begin{equation}
   \label{resG}
P_m =  \P_m
 \end{equation}  by induction on $m$. 

The $m=3$ case was checked in the introduction.

It will be convenient to use the following notation. 
Consider an $m$-dimensional vector space $\Rm$ with a  
basis $v_1,\ldots,v_m$. Then we can identify the variables $ y_{ijk}$
in  $P_m$  and  $\P_m$ with the standard 
basis of the space of three-vectors
 $$ y_{ijk} = v_i \wedge v_j\wedge v_k \in  \Lambda^3\Rm$$ 
and consider both polynomials $P_m$ and  $\P_m$
as elements of the space
$$  S^{({m-1})/{2}}(\Lambda^3\Rm)~.$$

There is a natural embedding of this space into the tensor algebra of $\Rm$,
 {\em i.e.,} the free associative algebra generated by $v_1, \dots, v_m$. 
Therefore  
we may view both
 the   left- and the right-hand  side of~(\ref{resG}) 
as polynomials in the \emph{non-commuting} variables $v_i$ 
\begin{align*}
P(\Gamma_m,\o_{can})&=P_m( y_{ijk})=P_m(v_1,v_2,\ldots, v_m)~,\\ 
\Pf(\Lambda(\Gamma_m)^{(1)}) &=\P_m( y_{ijk})=\P_m(v_1,v_2,\ldots,v_m)~.  
\end{align*}

The proof that  $P_m=\P_m$ proceeds  by showing that both $P_m$ and
$\P_m $ satisfy the same recursion formula. 

First, the generating function $P_m$. We have
\begin{equation}
  \label{eq:dc1}
P_m=P(\Gamma_m,\o_{can})=P_m'+P_m''~,\end{equation} 
where $P_m'$ generates the spanning trees of
 $\Gamma_m$ which contain the edge $\{1,2,3\}$, and $P_m''$ generates
those which do not. Given a spanning tree of the first kind, we can
collapse the edge $\{1,2,3\}$ to a new vertex, $0$, say, and  
obtain a spanning tree in
the complete $3$-graph $\Gamma'$ with $m-2$ vertices
$0,4,5,\ldots, m$.
 Conversely, every spanning tree in $\Gamma'$
can be lifted (in a non-unique way) to $\Gamma_m$, and every such lift  
together with the edge $\{1,2,3\}$ constitutes a spanning tree of
$\Gamma_m$.   This correspondence gives a
 relation between the generating functions which 
can be conveniently expressed in terms of the variables $v_i$:  
\begin{equation}
  \label{eq:dc2}
P_m'= y_{123} \, P_{m-2}(v_1+v_2+v_3, v_4, \ldots, v_m)~.  
\end{equation}
On the other hand, 
notice that  $P_m''$ is just the polynomial obtained from $P_m
( y_{ijk})$ by  setting $ y_{123}=0$.  

Now for the Pfaffian polynomial $\P_m$. 
Recall that $\Lambda(\Gamma_m)$ is the $m\times m$ matrix 
$$
\Lambda=(\lambda_{ij}), \quad \text{with} \quad \lambda_{ij}=\sum_k
 y_{ijk}, \ 1\leq i,j\leq m~.
$$  
 Write 
\begin{equation}
  \label{eq:dc3}
\P_m=\Pf(\Lambda^{(1)})=\P_m'+\P_m''~,
  \end{equation} 
where  $\P_m''$ is, \emph{by definition,} the 
polynomial obtained from $\P_m$ by 
setting $ y_{123}=0$.  We claim that   
\begin{equation} \label{eqR}
\P'_m= y_{123}\,\P_{m-2}(v_1+v_2+v_3, v_4, \ldots, v_m)~.
\end{equation} 
To show this, we apply the row development formula~(\ref{pf6})
to $\Lambda^{(1)}$ to obtain 
\begin{align*}
 \P_m&=\Pf(\Lambda^{(1)})\\
&=\lambda_{23}\,\Pf(\Lambda^{(1,2,3)}) -\lambda_{24}\,\Pf(\Lambda^{(1,2,4)})
+\ldots +\lambda_{2,m}\,\Pf(\Lambda^{(1,2,m)})~.
\end{align*} 
The only entry in
$\Lambda^{(1)}$ where $ y_{123}$ appears 
is 
$\lambda_{23}=-\lambda_{32}$, so that setting
$ y_{123}=0$ affects only the very first term of this expansion. 
It follows that 
$$\P'_m=\P_m-\P_m''\\
= y_{123}\,\Pf(\Lambda^{(1,2,3)})~.
$$
It remains to show that $\Pf(\Lambda^{(1,2,3)})$ 
is equal to 
$\P_{m-2}(v_1+v_2+v_3,v_4, \ldots, v_m)$. 
To see this, consider again the complete $3$-graph $\Gamma'$
with $m-2$ vertices $0,4,5, \dots , m$. 
The $(i,j)$ entry of the associated matrix is  
$$ y_{ij0}+\sum_{k=4}^m  y_{ijk}~.$$ 
Thus, if we delete the $0$-th row and column and 
substitute $v_0=v_1+v_2+v_3$, we get exactly the matrix 
$\Lambda^{(1,2,3)}$ (whose $(i,j)$ entry is $ y_{ij1}+ y_{ij2}+
y_{ij3} +\sum_{k=4}^m  y_{ijk}$). This proves~(\ref{eqR}).
\medskip

Thus we have shown that both $P_m$ and $\P_m $ satisfy the same
recursion relation
\begin{equation}
  \label{eq:recursion}
  P_m= y_{123}P_{m-2}(v_1+v_2+v_3,v_4,\ldots,v_m) 
      + \left[P_m\right]_{ y_{123}=0} 
\end{equation}
(and similarly for $\P_m$.)

This implies that $P_{m}=\P_{m}$ as follows. Since
$P_{m-2}=\P_{m-2}$ by the induction hypothesis, the
recursion~\eqref{eq:recursion} 
shows that $P_m - \P_m $ is divisible by $ y_{123}$. 
Since a similar recursion obviously holds 
for every edge  $\{i,j,k\}$, the difference
 $P_m - \P_m $ must be divisible by every $ y_{ijk}$.
Therefore,  $P_m - \P_m $ must be zero  
by degree count.  
This
completes the proof. \eeproof   
\medskip

\section{Properties of the Pfaffian-tree polynomial
$\P_m$}\label{sec:prop}     

In this section, we establish some algebraic properties of the
Pfaffian-tree polynomial of the complete $3$-graph $\Gamma_m$.

\subsection{Antisymmetry} 
\

By definition, $\P_m$ is a homogeneous
polynomial of degree $(m-1)/2$ in the indeterminates $
y_{ijk}$. Thinking of the $ y_{ijk}$ as elements of 
$\Lambda^3 \Rm$,
as in the proof of Theorem~\ref{mainth}, we may consider 
$\P_m$ as an element of the space 
 $ S^{(m-1)/2}\Lambda^3 \Rm.$ 
The following result shows that $\P_m$ belongs to the subspace
   $$(S^{(m-1)/2}\Lambda^3 \Rm)^-~,$$ 
where the superscript ${}^{-}$ indicates
 the subspace 
 which is totally antisymmetric  with respect to the
 action of the symmetric group $\Sy_m$. 
Recall that $v_1\ldots,v_m$ denotes a basis of $\Rm$ such that 
$ y_{ijk}=v_i\wedge v_j\wedge v_k$. 

\begin{proposition}\label{prop:alternate} 
For every permutation $\s\in\Sy_m$, one has  
  \begin{equation}
    \label{eq:alt}
\P_m(v_{\s(1)},\ldots, v_{\s(m)})= (-1)^\s \P_m(v_{1},\ldots, v_m).
    \end{equation}
\end{proposition}

\proof It is enough to prove this when $\s$ is a transposition. Thus,
we must show that $\P_m$ changes sign if two entries $v_i$ and $v_j$
are permuted.  In terms of the definition of $\P_m$ as a Pfaffian
(see~\eqref{eq:pf-tree-pol}), this follows from the fact 
that the Pfaffian of 
a skew-symmetric 
matrix changes sign if one
simultaneously permutes the $i$th and $j$th row and the $i$th and
$j$th column (see~\eqref{pf5}). 
Alternatively, another proof can be given if $\P_m$ is
viewed as the spanning tree generating function
$P(\Gamma_m,\o_{can})$. There  the proof comes down to the fact that
permuting two vertices reverses the orientation of $\Gamma_m$. We
leave the details of this alternative argument to the reader. \eeproof 

\subsection{Contraction-deletion relation}
\label{sec:contr-del}

\

Recall the recursion formula~\eqref{eq:recursion} shown in the proof
of Theorem \ref{mainth}:  
\begin{equation*}
  \label{eq:recursion2}
  \P_m= y_{123}\P_{m-2}(v_1+v_2+v_3,v_4,\ldots,v_m) +
  \left[\P_m\right]_{ y_{123}=0} 
\end{equation*} This formula can be viewed as a
  \emph{contraction-deletion} relation, as we shall now explain.   
It can be written as 
\begin{equation*}  \label{eq:dc_compl}
P(\Gamma_m,\o_{can}) = y_{\tilde e}P(\Gamma_m/e,\o_{can}/\o_{\tilde
 e})  +P(\Gamma_m - e,\o_{can})~, 
\end{equation*} 
where the
 notation is 
as follows. We have denoted by $e$   the edge
$\{1,2,3\}\in E(\Gamma_m)$, and by $\o_{\tilde e}$ the orientation of
the oriented edge $\tilde e =(123)$. The notation  $\Gamma_m - e$
stands for the $3$-graph $\Gamma_m$ with  
the edge $e$ deleted, and $\Gamma_m / e$ is the $3$-graph obtained
from $\Gamma_m$  by contracting the subgraph induced by the edge $e$
(that is, by replacing the three vertices of the edge $e$ in
$V(\Gamma_m)$ by a new vertex, say, $0$,  and discarding all edges
from $E(\Gamma_m)$ that become degenerate after this
identification). Notice that the quotient $3$-graph  
$\Gamma_m/e$ has multiple edges and, therefore, is not isomorphic to 
$\Gamma_{m-2}$.  Finally, $\o_{can}/\o_{\tilde e}$ is the induced
orientation of $\Gamma_m / e$. In our example, it is represented by
the cyclic permutation $(045\ldots m)$ of the vertex set
$V(\Gamma_m/e)$.  

This \emph{contraction-deletion} relation can be formulated in
general: 
\begin{proposition}[{\bf Contraction-deletion relation}] 
\label{prop:contr-del} 
\begin{equation}
  \label{eq:3del-cont}
 P(G,\o) =  y_{\tilde{e}} P(G/e, \o/\o_{\tilde e})+P(G - e, \o)~,   
\end{equation}
where $G$ is a $3$-graph with an orientation $\o$ and an oriented edge  
$\tilde e= (e,\o_{\tilde e})$ and $\o/\o_{\tilde e}$ is the induced
orientation of the quotient $3$-graph $G/e$. 
\end{proposition}

We omit the easy proof, and merely spell out the rule to compute the
induced orientation $\o/\o_{\tilde e}$. Assume
$V(G)=\{1,2,\ldots,m\}$,  and, as above, let $0$ be the new vertex of
$G/e$ obtained by contracting the vertices of $e$. 
Then $\o/\o_{\tilde e}$ is represented by a cyclic permutation
$\sigma$ of $V(G/e)$ such that if one inserts a cyclic permutation
representing $\o_{\tilde e}$ 
in place of $0$ into $\s$, one obtains a cyclic permutation
representing the original orientation $\o$ of $G$. For example,
$$
(12345\ldots m)/(124)=-(12435\ldots m)/(124)=-(035\ldots m)
$$ 
which shows that the analogue of (\ref{eq:recursion}) for the edge
$\{1,2,4\}$ is 
\begin{equation}   \label{eq:recursion3}
\P_m=- y_{124}\P_{m-2}(v_1+v_2+v_4,v_3,v_5,\ldots,v_m) +
\left[\P_m\right]_{ y_{124}=0}
\end{equation}

\begin{remark}
{\em A contraction-deletion relation analogous to~(\ref{eq:3del-cont}) 
exists for ordinary graphs and relates the 
Kirchhoff polynomial $\mathcal{D}_G$ of a graph $G$ to those of $G/e$
and $G-e$, the graphs  obtained from $G$ by, respectively, 
contracting and deleting an edge $e\in E(G)$. In fact, one of the
standard proofs of the classical Matrix-Tree Theorem is based on this
relation (see {\em e.g.,} the proof given in \cite[Theorem
II.12]{Bol}). Note that our proof of Theorem~\ref{mainth} is similar
in spirit. The Pfaffian-tree polynomial $\P_G$ satisfies a
contraction-deletion relation analogous to~\eqref{eq:3del-cont}, and a
possible variant of our proof would be to prove this relation for
$\P_G$ independently and then deduce Theorem~\ref{mainth} from it.
However, we have found it more convenient to derive
Theorem~\ref{mainth} from its particular case of complete $3$-graphs:
in our proof, the contraction-deletion relation for complete
$3$-graphs takes the form of equations \eqref{eq:dc1} and
\eqref{eq:dc2} for the tree generating function, and of
equations~\eqref{eq:dc3} and~\eqref{eqR} for the Pfaffian,
respectively. } 
\end{remark}

\subsection{Three-term relation}

 \begin{proposition} \label{three-term} The polynomial $\P_m$
 satisfies the relation  
 \begin{equation}  \label{eq3t}
 \P_m(v_2+v_3,v_4,\ldots) +\P_m(v_3+v_4,v_2,\ldots)
 +\P_m(v_2+v_4,v_3,\ldots) =0 
 \end{equation}
  where the dots stand for $v_5,v_6,\ldots, v_{m+2}$.
 \end{proposition}

 \proof Here we think of $\P_m$ as the tree-generating function 
$P(\Gamma_m,\o_{can})$. The first summand in (\ref{eq3t}) is obtained
from $\P_m(v_2,v_4,\ldots)$ by replacing every $ y_{2ij}$ occurring in
it by $ y_{2ij} + y_{3ij}$ and expanding by multilinearity. In
particular, a monomial in $\P_m(v_2,v_4,\ldots)$  corresponding to a
tree $T$ gets replaced by $2^{n}$ terms, where $n$ is the valency of
the vertex $2$ in $T$. If the other two summands of (\ref{eq3t}) are
expanded similarly, then each of these $2^{n}$ terms coming from the
tree $T$ also occurs in exactly one of the two other summands, but
with opposite sign (by the 
antisymmetry
of $\P_m$, see (\ref{eq:alt})). Thus, all terms in (\ref{eq3t}) cancel, 
and the result follows.  
\eeproof

The three-term relation implies the following properties of $\P_m^2$
which play an important role in our study of the lowest order term of
the Alexander-Conway polynomial in  \cite{MV2}.  

\begin{corollary} \label{cor.recursion}   The polynomial $\P_m^2$
satisfies the following relations: 
$$\begin{array}{lcl}
(i)\hspace{ -8pt} &\displaystyle{\left[\frac {\partial^2 \,
\P_{m+2}^2}{\partial  y_{123}^2}\right]_{v_1=0}} 
&\hspace{ -15pt}= 2\, \P^2_{m}( v_2+v_3,\ldots)\\
&\vspace{-10pt}&\\
(ii)\hspace{ -8pt} &\displaystyle{\left[\frac {\partial^2 \,
\P_{m+2}^2}{\partial  y_{123}\, \partial  y_{124}}\right]_{v_1=0}} 
&\hspace{ -15pt}= \P^2_{m}( v_2+ v_3, v_4, \ldots) +\P^2_{m}
( v_2+v_4,v_3, \ldots) \\ 
&& \hspace{ 20pt} - \P^2_{m}( v_3+v_4,v_2, \ldots) \\ 
(iii)\hspace{ -8pt} &\displaystyle{\left[\frac {\partial^2 \,
\P^2_{m+2}}{\partial  y_{123}\, \partial  y_{145}}\right]_{v_1=0}}
&\hspace{ -15pt}= \P^2_{m}( v_3+v_4,v_2, v_5, \ldots) + \P^2_{m}( v_2
+v_5, v_3, v_4, \ldots)\\ 
&&\hspace{ -10pt}- \P^2_{m}( v_2+v_4,v_3, v_5, \ldots) - \P^2_{m}( v_3
+v_5, v_2, v_4, \ldots) 
\end{array}
$$ 
where the dots stand for the $v_i$ with indices not
involved on the left-hand side (for example,  in the first equation, 
the dots stand for $v_4,v_5, \ldots, v_{m+2}$.)
 It also satisfies
all equations obtained from the above ones by some permutation of 
the indices $1,2,\ldots, m+2$.  
\end{corollary}
\proof 
The contraction-deletion formula~\eqref{eq:recursion} shows
that 
$$\frac {\partial \, \P_{m+2}}{\partial  y_{123}}
=\P_{m}(v_1+v_2+v_3,v_4,\ldots).$$ 
Since
$\left[\P_{m+2}\right]_{v_1=0}=0$, it follows that 
$$\left[\frac
{\partial^2 \, \P_{m+2}^2}{\partial
y_{123}^2}\right]_{v_1=0}=2\,\left( 
\left[\frac {\partial \, \P_{m+2}}{\partial
y_{123}}\right]_{v_1=0}\right)^2= 2\, \P^2_{m}( v_2+v_3,v_4,\ldots)~,$$
proving relation (i). 

For relation (ii), we have 
\begin{align*} 
\left[\frac {\partial^2 \, \P_{m+2}^2}{\partial
y_{123}\, \partial  y_{124}}\right]_{v_1=0}&= 2\,\left[\frac {\partial
\, \P_{m+2}}{\partial  y_{123}}\right]_{v_1=0} \left[\frac {\partial
\, \P_{m+2}}{\partial  y_{124}}\right]_{v_1=0} \\ 
&= - 2\,\P_{m}(v_2+v_3,v_4,\ldots)\P_{m}(v_2+v_4,v_3,\ldots)~,
\end{align*} 
where the minus sign comes from the minus sign in
\eqref{eq:recursion3}. Since
$$
\P_{m}(v_2+v_3,v_4,\ldots)+\P_{m}(v_2+v_4,v_3,\ldots) =
-\P_{m}(v_3+v_4,v_2,\ldots)
$$   
by the three-term relation~(\ref{eq3t}), relation (ii) now follows  
from the identity 
\begin{equation}
  \label{eq:iden}
 -2AB=A^2+B^2-(A+B)^2. 
\end{equation}

The proof of relation (iii) is similar in spirit but more
complicated. Let us abbreviate $\P_{m}(v_2+v_3,v_4,v_5,\ldots)$ by
$\P(2+3,4,5)$, and similarly for the other terms. After computing
derivatives as above, the left-hand side of relation~(iii) is
\begin{equation}\label{eq:lhs}
2\,\P(2+3,4,5)\P(4+5,2,3)~.
\end{equation}
Applying the three-term relation to
the second factor, and using the 
antisymmetry 
of the polynomial $\P_m$ (see~\eqref{eq:alt}), 
we see that~(\ref{eq:lhs}) is equal to
\begin{align*} 
-2\,&\P(2+3,4,5)\bigl(\P(5+2,4,3) + \P(2+4,5,3)\bigr) \\
&=-2\,\P(2+3,5,4)\P(2+5,3,4) + 2\P(2+3,4,5)\P(2+4,3,5)~.
\end{align*}
Applying the three-term relation again to the two products above, 
and using~(\ref{eq:iden})  as in the proof of~(ii), 
we find that~(\ref{eq:lhs}) is equal to  
\begin{align*} & \bigl( \P(2+3,5,4)^2 +\P(2+5,3,4)^2 -
\P(3+5,2,4)^2\bigr) \\  
&\ \ \ \ \ - \bigl(\P(2+3,4,5)^2+\P(2+4,3,5)^2 - 
\P(3+4,2,5)^2\bigr)\\ 
=&\ \P(2+5,3,4)^2+\P(3+4,2,5)^2-\P(2+4,3,5)^2-\P(3+5,2,4)^2
\end{align*}
which is the right-hand side of relation~(iii).  This completes the
proof. 
\eeproof

Relations (i)-(iii) can be used to compute $\P_m^2$ recursively. This
relies on the following simple observation.  
\begin{lemma} \label{twice}
  If a monomial 
$\prod_\alpha  y_{i_\alpha j_\alpha k_\alpha}$
occurs with non-zero coefficient in $\P^2_m$, then there exists
$p\in\{1,2,\ldots, m\}$ such that $p$ occurs exactly twice in the list
of indices $i_1,j_1,k_1,i_2,\ldots j_{m-1},k_{m-1}$. 
\end{lemma}
\proof 
Recall that  $\P_m$ is the tree generating function $\sum_T
y(T,\o_{can})$ of the complete $3$-graph $\Gamma_m$. A monomial
$M=\prod_\alpha  y_{i_\alpha j_\alpha k_\alpha}$ determines a
$3$-graph $G_M$ whose edges are given by the $ y_{ijk}$ occurring in
$M$.  The  coefficient of $M$ in $\P_m^2$ is the number (with signs)
of ways  the $3$-graph $G_M$ can be written as the union of two
spanning trees $T$ and $T'$. Since every vertex  is incident with at
least one edge of $T$ and one edge of $T'$, every  $p\in\{1,2,\ldots,
m\}$ occurs at least twice in the list of indices
$i_1,j_1,k_1,i_2,\ldots j_{m-1},k_{m-1}$ of $M$. Since the total
number of indices in this list is $<3m$, there must be an index which
occurs exactly twice. 
\eeproof 

\begin{corollary}  Relations (i)-(iii) of Corollary~\ref{cor.recursion}
together with Lem\-ma~\ref{twice} allow to compute 
$\P_m^2$ recursively with initial condition $\P_1=1$. 
  \end{corollary}
\proof For every monomial occurring with non-zero coefficient in
$\P_m^2$, we can find by Lemma \ref{twice} an index $p$ such that the
monomial contains $ y_{pij} y_{pkl}$ for some $i,j,k,l$, but no other
$ y_{\alpha\beta\gamma}$ with $p\in\{\alpha,\beta,\gamma\}$. The
coefficient of such a monomial in $\P_m^2$ 
is equal to its coefficient in 
$$ y_{pij}\, y_{pkl}\, \left[\frac {\partial^2 \, \P_{m}^2}{\partial
     y_{pij}\, \partial  y_{pkl}}\right]_{v_p=0}$$ 
if $\{i,j\}\neq\{k,l\}$,
  or to one half of this coefficient if $\{i,j\}=\{k,l\}$. 
But this coefficient can be computed
recursively by relations (i)-(iii) of
Corollary~\ref{cor.recursion}.\nolinebreak\eeproof   

\begin{remark}\rm 
Monomials in $\P_m$ correspond to spanning trees and always occur with 
coefficient $\pm1$. 
Therefore the coefficient of a monomial $M$ in  $\P_m^2$ is equal to 
the  number (with signs) of ordered tree decompositions of the 
associated $3$-graph $G_M$, divided by the symmetry factor $|Aut(G_M)|$.
Here,  by an ordered  tree decomposition of $G_M$ we mean
a sub-$3$-graph $T$ which is a tree and whose complement is also a tree, 
and by $Aut(G_M)$ the group of automorphisms of $G_M$
inducing the identity map on the set of vertices $V(G_M)$. 
The cardinality  $|Aut(G_M)|$  is equal to $2^{d}$, 
where  $d$ is the number of (unordered) 
triples of vertices in $G_M$ 
with $2$ edges  attached to them.
(This is because there can be at most two edges 
attached to a triple of vertices of  $G_M$ if $M$ has non-zero 
coefficient in $\P_m^2$.) 

Here are three examples to illustrate this.   
The monomial $y_{123}^2$ has two ordered tree decompositions and 
$|Aut(G_M)|=2$; therefore it occurs with coefficient $1$ in $\P_3^2$. The
monomial $y_{123}^2y_{245}y_{345}$ has four ordered tree decompositions (each
occurring with a plus sign) and $|Aut(G_M)|=2$; therefore it occurs with
coefficient $2$ in $\P_3^2$. Finally, the monomial  
$$ M=y_{145}\, y_{146}\, y_{256}\, y_{257}\, y_{347}\, y_{367}$$ 
has six  ordered tree decompositions (again each contributing $+1$) and
$|Aut(G_M)|=1$; therefore it occurs with coefficient $6$ in $\P_7^2$. 
\end{remark}
 
\

\subsection{Four-term relation}\label{sec:four-term-rel}

\begin{proposition} 
\label{four-term} 
The polynomial $\P_m$ satisfies the relation 
\begin{equation}  \label{murel}
{\partial \P_m \over \partial  y_{ijk}} -{\partial \P_m
    \over \partial  y_{ijl}}={\partial \P_m \over
    \partial  y_{jkl}}-{\partial \P_m \over \partial y_{ikl}} 
\end{equation}
 for every set $\{i,j,k,l\}$ of four distinct vertices.
\end{proposition}

See Figure~\ref{figmu} for a pictorial illustration of~(\ref{murel}).

\begin{figure}[h]   
\begin{center}
\MUa\  $-$ \ \MUb \ = \ \MUc \ $-$ \ \MUd

\caption{\label{figmu} An illustration of Equation~(\ref{murel}).}
\end{center}
\end{figure}

\medskip

\proof The four-term relation \eqref{murel} can be deduced from the
three-term relation, as follows.  

Without loss of generality, we may assume that $(i,j,k,l)=(1,2,3,4)$. 
By the contraction-deletion relation~(\ref{eq:3del-cont}), 
the left-hand side  of~(\ref{murel}) is 
$$
\P_{m-2}(v_1+v_2+v_3,v_4,\ldots) + \P_{m-2}(v_1+v_2+v_4,v_3,\ldots)
$$
and the right-hand side is equal to 
$$-\P_{m-2}(v_2+v_3+v_4,v_1,\ldots) - \P_{m-2}(v_1+v_3+v_4,v_2,\ldots)~.$$
(The signs 
are obtained as in \eqref{eq:recursion3};  see
the discussion following Proposition \ref{prop:contr-del}.)
Using the three-term relation~(\ref{eq3t}) 
and antisymmetry, we see that both sides are equal to 
$$\P_{m-2}(v_1+v_2,v_3+v_4, \ldots)$$
which proves~(\ref{murel}). \eeproof

Here is an  equivalent formulation of Proposition~\ref{four-term},
which we use in~\cite{MV2}. 

Let us consider the polynomial $\P_m\in S^{(m-1)/2} \Lambda^3 \Rm$ as
a {\em polynomial  function} of degree $(m-1)/2$ on the vector space
$W=\Lambda^3 V$, where $V$ is the dual of $\Rm$. 
Note that  $\{y_{ijk}\}$ is a basis of the space of
linear  forms on $W$. Let $\Y_{ijk}\in W$
be the dual basis, {\em i.e.\/} the evaluation
$$
 \langle y_{ijk},\Y_{\alpha\beta\gamma}\rangle
$$ 
is the sign of the permutation
$({}^{i}_\alpha{}^j_\beta{}^k_\gamma)$ if
$\{i,j,k\}=\{\alpha,\beta,\gamma\}$, and is zero otherwise.   

\begin{proposition} \label{8.2} 
The polynomial $\P_m$
descends to a well-defined polynomial function on $W/W_0$, 
where $W_0$ is the subspace of $W$ generated by vectors of the form   
$$
(\Y_{ijk} -\Y_{ijl}) - (\Y_{jkl}-\Y_{ikl})
$$ 
for every set $\{i,j,k,l\}$ of four distinct vertices.  
\end{proposition}
 \noindent
\proof
 Equation~(\ref{murel})   in~\ref{four-term} shows that the derivative 
 of $\P_m$ in the direction of any vector in $W_0$  is identically
 zero. This implies Proposition~\ref{8.2} by Taylor's
 formula. \eeproof  
\medskip

\noindent
\begin{remark} {\rm
The four-term relation for the 
tree-generating function 
$\P_m=P(\Gamma_m,\o_{can})$ has a simple combinatorial meaning, as 
follows. 

The
partial derivative ${\partial \P_m / \partial  y_{123}}$  
       is equal to $ y_{123}^{-1} \P'_m,$ where $\P'_m$ is
the generating function for those trees which contain the edge
$\{1,2,3\}$. 
Given such a tree $T$, let $T'$ denote $T$ with the edge $\{1,2,3\}$
removed. Note that  $T'$ is the disjoint union of $3$ subtrees 
$$
T'=T_1\cup T_2\cup T_3~,$$ 
where $T_i$ 
denotes the subtree containing 
the vertex $i$.  The key point is to observe that 
gluing the edge $\{1,2,4\}$ (resp.\ $\{1,3,4\}$, \
$\{2,3,4\}$) to $T'$
yields a tree if and only if the vertex $4$ is 
contained in the component $T_3$ (resp.\ $T_2$, \
$T_1$).

It follows that for every tree  contributing to 
${\partial \P_m / \partial y_{123}}$, there is a unique way 
of replacing the edge $\{1,2,3\}$ by one of the three other edges,  
    $\{1,2,4\}$  $\{1,3,4\}$, or $\{2,3,4\}$, so that the result is
again a  tree. 
Note that the tree thus obtained contributes to one and only one of 
 ${\partial \P_m / \partial  y_{124}}$,  ${\partial \P_m / \partial
  y_{134}}$, and  ${\partial \P_m / \partial  y_{234}}$.  

This observation already  implies~(\ref{murel}) \emph{modulo} $2$. 
Indeed,  
the trees contributing to each of the four partial derivatives 
${\partial \P_m / \partial  y_{ijk}}$ are partitioned
into three
disjoint subsets;
the set of these altogether $12$ subsets is divided 
into six pairs of two; the two subsets in every pair are in bijective
correspondence with each other (so, in particular, they have the same
cardinality), but contribute to different 
${\partial \P_m / \partial y_{ijk}}$. This is the combinatorial
meaning of the  four-term relation up to sign.   

The signs can also be checked combinatorially, as follows.
Let us do for example the
case corresponding to  the pair 
$\{(123),(124)\}$. In this case,  we must look at a tree $T$ which is
the union of the edge $\{1,2,3\}$ and a remainder, $T'$, so that $T'$
plus the edge $\{1,2,4\}$ gives again a tree, $ \widehat T$,
say. (This means that the vertex $4$ must lie in the component $T_3$
of $T'$.) Thus, we have   
\begin{equation}
{\partial y(T, \o_{can})\over \partial y_{123}}\ = \ \pm \  {\partial 
  y(\widehat T , \o_{can})\over \partial y_{124}} \label{sig} 
\end{equation}
  and we must show that the sign is $+1$ in this case. This can be
seen
by comparing the orientations induced by planar embeddings
  of $T$ and $\widehat T$ which coincide on $T'$.

\begin{figure}[h]
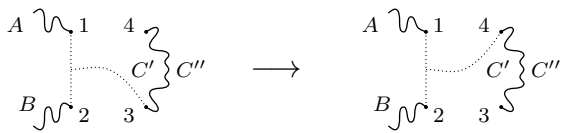

\begin{center}
\FIGPFi  
$\longrightarrow$ 
\ \ \ \ \ \ \ \  \ \ \FIGPFii

\caption{\label{figpf} Replacing $(123)$ with $(124)$ in a planar
embedding  of $T$ gives a planar embedding of $\widehat T$. 
} 
\end{center}
\end{figure}

We proceed as
in the proof of the sign change property in 
Proposition \ref{signs}. From the embedding of $T$ we may read 
off an ordering $\o$ of the form $1A2B3C$ (see the left diagram  
on Figure~\ref{figlem}). But now we need to take into account where
the vertex $4$ is placed. Therefore we decompose $C$ as $C''4C'$ (see 
Figure~\ref{figpf}) and write $\o=1A2B3C''4C'$. The corresponding
embedding of $\widehat T$ gives an ordering $\hat\o$ of the form 
$1A2B4C'3C''$. Observing that $C'\cup C''$ has an odd number of
elements, it is easy to check that $\o$ and $\hat\o$ induce the same
orientation, showing that the sign in~(\ref{sig}) is indeed $+1$, as
asserted.  

The recursion formula for $\P_m^2$ of Corollary \ref{cor.recursion}
can be understood combinatorially in a similar way.
}
\end{remark}

\end{document}